\input AHTOHFIE.STY
\hfuzz 11pt
\UDC{
512.543.52 +  
512.543.72    
}
\def\NN{{\cal N}}
\def\M{{\cal M}}
\def\A{{\cal A}}
\def\B{{\cal B}}
\def\fp{{\bf fp}}

\title{
Commutator length of powers in free products of groups
}

\author{%
Vadim Yu. Bereznyuk
\quad
\quad
Anton A. Klyachko
}
\address{
Faculty of mechanics and mathematics of Moscow State University
\\
Moscow 119991, Leninskie gory, MSU.
\\
Moscow Center for Fundamental and Applied Mathematics.
\\
kuynzereb@gmail.com
\quad
\quad
klyachko@mech.math.msu.su
}
\grants{\CENTR075-15-2019-1621
and also
\RFBR19-01-00591}

\abstract{\narrower
Given groups $A$ and $B$,
what is the minimal commutator length
of the
\the\year th (for instance) power of an element~$g\in A*B$ not conjugate
to elements of the free factors?  The exhaustive answer to this question
is still unknown, but we can give an almost answer:  this minimum is one
of two numbers (simply depending on $A$ and $B$). Other similar problems
are also considered.
}

\s 0.
Introduction

It is well known that, in free groups, nonidentity commutators are not
proper
powers~[Sch59].
A product of two commutators in a free group can surely be
the square of a nonidentity element, and can even be a cube,
as Culler noticed {[Cull81]}:
$
[a,b]^3 = [a^{-1} ba, a^{-2}bab^{-1}][bab^{-1}, b^2].
$
This equality holds in the free group $F(a,b)$ and, therefore,
for any elements $a$ and $b$ of any group.
Moreover, Culler {[Cull81]} showed that, in the free group $F(a,b)$,
the element
$[a,b]^n$
decomposes into a
product of $k$ commutators if $n\le 2k-1$.

For free groups,
Culler's estimate cannot be improved in any sense:
\disp{\sl\narrower\narrower\narrower\narrower
if, for some elements $x_i,y_i,z$ of a free group,
$[x_1,y_1]\ldots[x_k,y_k]=z^n$,
where $n\ge2k$, then $z=1$.
}%
This remarkable fact was obtained in {[CCE91]} for $k=2$
and in {[DH91]} in the general case.
In the same paper~{[DH91]}, a similar assertion
was proven
for free products of \emph{locally indicable}
groups (i.e. groups, in which each nontrivial finitely
generated
subgroup admits an epimorphism onto $\Z$). Later,
it was discovered that this assertion remains valid
in free products of any torsion-free groups:
\disp{\sl
if\/ $[x_1,y_1]\ldots[x_k,y_k]=z^n$ for some elements $x_i,y_i,z$ of a free
product of torsion-free groups,
where $n\ge2k$, then $z$ is conjugate to an element of a free factor.
}%
This was shown in [Ch18] and [IK18] (independently).
Moreover, both papers
mentioned that the arguments remain valid
if the torsion-free condition is replaced with a
small-torsion-free condition. However, arguments in
[Ch18] and [IK18] are different:
\-
Chen's proof is based on Calegari's approach [Cal09],
\-
while [IK18]
is based on the car-crash lemma [K93];

\enditem
this is why, the results of [Ch18] and [IK18]
for groups with torsion
are
different
(and even incomparable --- neither one is
stronger than the other):
\disp{\sl
suppose that, for some elements $x_i,y_i,z$
of a free
product of groups
without nonidentity elements of order less than
$N$,
an equality
$[x_1,y_1]\ldots[x_k,y_k]=z^n$
holds\;
then $z$ is conjugate to an element of a free factor
}
\vskip-9mm
$$
\sl
\qqbox{if}
\cases{
n
\ge
2k+\[\frac{2n}N\]&\rm[Ch18]
\qbox{\(henceforth, $[x]\:=\max\{p\in\Z\;|\;p\le x\}$\)}
\cr
\hbox{\sl or}&
\cr
n\ge2k\qqbox{and}
N>n&\rm[IK18].
}
\eqno{(*)}
$$

We show that condition $(*)$ can be replaced with
a weaker condition
$$
n
\ge
2k+2\[\frac{n}N\].
\eqno{(**)}
$$
Clearly, this strengthens both results $(*)$.
Moreover,
estimate $(**)$ is the best possible or almost the best possible.
Namely, the situation is as follows.

Let $G$ be a group
with a fixed free-product decomposition:
$G=\zvezda\limits_{j\in J} A_j$.
Let
$k(G,n)$ be the minimal $k\in\Z$ such that the $n$th
power of an element of $G$ not
conjugate to elements of~$\bigcup\limits_{j\in J} A_j$
decomposes into a product of $k$
commutators and let $N(G)$ be the minimal order of a nonidentity
element of~$G$. Thus, according to $(**)$,
any free product $G$
satisfies the inequality
$k(G,n)\ge\[{n\over2}\]-\[{n\over N(G)}\]+1$.
This estimate is almost the best possible in the following sense:
Theorem~1 (see the following section)
asserts, in particular, that
\disp{\sl
for any free product $G=\zvezda\limits_{j\in J} A_j$,
the value $k(G,n)$ is
$
\hbox{\sl either}\quad
\[{n\over2}\]-\[{n\over N(G)}\]+1
\qqbox{\sl or}
\[{n\over2}\]-\[{n\over N(G)}\]+2.
$
}
Putting $k(G,n)=1$, we obtain a well-known fact [CER94]:
\disp{\sl
in any free product $G$,
a commutator not conjugate to elements of the free factors
can be a proper power only if $N(G)=2$ or $N(G)=3$\;
in the latter case, this commutator can only be a cube.
}%
For larger $k(G,n)$, our result is (apparently) new.

\smallskip

\noindent
Actually, we study equations more
general than the equation $[y,z][t,u]\dots=x^n$ considered above:
\-
the power $x^n$ is replaced with a ``generalised power",
i.e. the product of conjugate elements;
\-
and the product of commutators is replaced with a product of
commutators and elements conjugate to elements of the free factors.

\proclaim Main theorem \rm(a simplified form).
Suppose that, in a
free product of groups $G=\zvezda\limits_{j\in J} A_j$
without
nonidentity elements of order less than $N$,
an equality
$$
c_1\dots c_kd_1\dots d_l=u_1^{n_1}\dots u_m^{n_m}
$$
holds,
where $c_i$ are commutators, $d_i$ are conjugate to elements of
$\bigcup\limits_{j\in J}A_j$,
elements
$u_i$ are conjugate to each other and not conjugate to elements of
$\bigcup\limits_{j\in J}A_j$,
and $n_i$ are positive integers. Then
$$
2k+l
\ge
\sum_{i=1}^m(n_i-1)
-2\[\frac1{N}\sum\limits_{i=1}^mn_i\]
+2.
$$

This result significantly strengthens earlier known facts:
\disp{\sl
under the hypothesis of the main theorem
$
2k+l\ge
\cases{
\sum\limits_{i=1}^m(n_i-1)
-\[\frac2{N}\sum\limits_{i=1}^mn_i\]+2,
&if $l=0$ \quad\quad\quad\rm[Ch18];
\cr\cr
\sum\limits_{i=1}^m(n_i-1)
+2,
&if $N>\sum\limits_{i=1}^mn_i$\quad\rm[IK18].
}
$
}%
The main theorem immediately implies what is said above on
inequality $(**)$.

\Corollary 1.
Suppose that, in a free product of groups $G=\zvezda\limits_{j\in J} A_j$,
an equality
$
c_1\dots c_k=u^n
$
holds,
where $c_i$ are commutators and $u$ is not conjugate to elements of
the free factors. Then
$
2k\ge n-2\[{n\over N}\]+1
$ {\rm(or,
equivalently,
$
k\ge [{n\over2}]-\[{n\over N}\]+1
$%
).}

The statement of the main theorem above is somewhat simplified.
In fact, we prove a stronger estimate under weaker
assumptions. The full statement of the main theorem and its proof
can be found in the last section.
In Section 2, we derive
Theorem~1 (mentioned
above) from the main theorem.
Sections 3 and 4 contain necessary information
about Howie diagrams and
motions on surfaces, i.e. about the car-crash lemma.
This lemma from~[K93] (or its variants) was already
applied in [FK12] and [IK18] to problems related to
the
commutator length
(and, e.g., in
[K93],
[ClG95],
[FeR96],
[Kl05],
[Kl06a],
[Kl06b],
[Kl07],
[Cl03],
[ClG01],
[CoR01],
[FoR05a],
[FoR05b],
[Kl09],
[Le09],
and
[KlL12],
it is applied to other problems).
We need a new version of the car-crash lemma, which is discussed in
Section 4. Surprisingly, a substantial role in that section
is played by the \emph{fair partition problem}, see, e.g., [Me06].

\s
Notation

{
Our notation}
is mainly standard. Note only that, if
$k\in \Z$, and $x$ and $y$ are elements
of a group, then $x^y$, $x^{ky}$,
and~$x^{-y}$ denote $y^{-1}xy$, $y^{-1}x^ky$, and $y^{-1}x^{-1}y$,
respectively.
The commutator~$[x,y]$ is $x^{-1}y^{-1}xy$.
The symbol~$\cl(g)$ denotes the
\emph{commutator length} of an element $g$ of a group,
i.e. $\cl(g)$ is the minimal integer $k$ such that
$g$ decomposes into a product of $k$ commutators
(and $\cl(1)=0$).
The word ``surface" always means a
closed surface (not necessarily connected).
The Euler characteristic of a surface $S$ is
denoted by~$\chi(S)$.
The letters $\R$, $\Z$, and $\N$ denote the set of
real, integer, and natural (positive integer) numbers,
respectively.
The symbol $[x]$ denotes the integer part of a real number~$x$
(i.e. $[x]$ is the maximal integer not exceeding $x$).

\s 1.
Powers of small
commutator length

Culler's bound mentioned in the very beginning of this paper
can be stated as follows.

\proclaim Culler's inequality {\rm[Cull81]}.
For any elements $a$ and $b$ of any group and
for any $n\in\N\cup\0$,
$$
\cl\([a,b]^n\)\le\[{n\over2}\]_c+1,
\qqbox{where}
[x]_c\:=\cases{
[x],&if $x\ne0$\;
\cr
-1, &if $x=0$.
}
$$

\Lemma 1.
If $a$ and $b$ are elements of a group
and $m\in\N$,
then $(ab)^m$ is conjugate to an element of the form
$
a^mb^mc_1c_2\dots c_{\[m\over2\]},
$
where $c_i\in G$ are commutators.

\Proof
$$
\eqalign{
&a^l(ba)^sb^l\cdot[a^{l-2}b^{l-1},b^{2-l}a^{1-l}]=
\cr
&=
a^l(ba)^sb^l\cdot
b^{1-l}a^{2-l}a^{l-1}b^{l-2}a^{l-2}b^{l-1}b^{2-l}a^{1-l}=
a^l(ba)^sbab^{l-2}a^{l-2}ba^{1-l}
\sim
a^{l-2}ba(ba)^sbab^{l-2}=a^{l-2}(ba)^{s+2}b^{l-2}.
}
$$
An obvious induction shows that, for some commutators
$c_i$, the element
$a^mb^mc_1c_2\dots c_{\[m\over2\]}$ is conjugate to
$(ba)^m$ if $m$ is even, or to
$a(ba)^{m-1}b$ if $m$ is odd. This completes the proof
(because $a(ba)^{m-1}b=(ab)^m\sim (ba)^m$, where $\sim$
means conjugation).

\Lemma 2.
If $a$ and $b$ are elements of a group,
$m\in\N\ni s$, and $a^m=b^m=1$,
then
$\cl\((ab)^{ms}\)\le s([m/2]-1)+\[{s/2}\]_c+1$.

\Proof
Note that, for any nonidentity element $g$ of
the commutator subgroup of
any group and for
any $s\in\N\cup\0$,
$$
\cl(g^s)\le s(\cl(g)-1)+\[{s\over2}\]_c+1.
$$
Indeed,
representing $g$ as $g=ch$, where $c$ is a
commutator, and $\cl(h)=\cl(g)-1$, we obtain
$$
\cl(g^s)=\cl((ch)^s)=\cl(c^sh^{c^{s-1}}h^{c^{s-2}}\dots h^ch)\le
\cl(c^s)+s\cdot\cl(h)\le
\[{s\over2}\]_c+1+s(\cl(g)-1)
\qbox{(by Culler's inequality)}.
$$
This completes the proof, because
$\cl\((ab)^m\)\le [m/2]$ by Lemma 1.

\Th 1.
For any free product $G=\zvezda\limits_{j\in J} A_j$
and any $n\in\N$,
$$
\eqalign{
&\hbox{either }
k(G,n)=\[{n\over2}\]-\[{n\over N(G)}\]+1
\qqbox{or}
k(G,n)=\[{n\over2}\]-\[{n\over N(G)}\]+2,
\qbox{where}
\cr
&k(G,n)\:=\min\left\{\cl(g^n)\;\Biggm|\;
g\in G,\
g \hbox{ is not conjugate to elements of }
\bigcup\limits_{j\in J} A_j\right\}
\hbox{ and }
N(G)\:=\min\bigl\{|\!\gp g\!|\;\bigm|\;g\in G\setminus\1\bigr\}.
}
$$
Moreover,
$k(G,n)=\[{n\over2}\]-\[{n\over N(G)}\]+1$ if
at least one of the following conditions is satisfied\:
\newline
{\rm a)}
$n$ is even and $\[{n\over N(G)}\]$ is odd\;
\qquad
\qquad
{\rm b)}
$n$ is divisible by $N(G)$\;
\qquad
\qquad
{\rm c)}
$n\le N(G)$;
\qquad
\qquad
{\rm d)}
$N(G)=2$.

\Proof
Surely, $N(G)$ is either prime or infinite.  For $N(G)=2$, the assertion
holds, because the group $G$ in this case contains an infinite dihedral
subgroup, whose commutator subgroup coincides with the set of commutators
(and trivially intersects conjugates of free factors). For
$N(G)=\infty$, the argument below is essentially valid, but we leave it to
readers, because the assertion of the theorem in this case follows
immediately from the results of [Ch18] and [IK18] mentioned in the
introduction.  Thus, we assume that $N(G)$ is odd.

If $z^m=1$, then we have two bounds:
$$
\cl([x,y]^m)
\le
\[{m\over2}\]_c+1
\qqbox{and}
\cl\([z,u]^{ms}\)
\le
s\(\[{m\over2}\]-1\)+\[{s\over2}\]_c+1,
\eqno{({**}*)}
$$
The first one is Culler's inequality,
and the second one is Lemma 2.

Consider in $G$ a commutator $[z,u]$, where $z^{N(G)}=1$ and
$u$ does not lie in the same free factors as $z$.
Let us divide $n$ by $N=N(G)$ with
remainder: $n=rN+t$, where $0\le t<N$
(and $r=\[{n\over N}\]$).
Let the symbols $\Delta_0(a,b,\dots)$ and $\Delta_{odd}(a,b,\dots)$
denote the
number of zeros and the number of odd numbers in the
tuple~$(a,b,\dots)$.
Then for odd $N$, we obtain
$$
\eqalign{
k(G,n)
&\le
\cl\([z,u]^n\)
=
\cl\([z,u]^{rN+t}\)
\le
\cl\([z,u]^{rN}\)+\cl\([z,u]^t\)
\lee^{{**}*}
\(r\(\[{N\over2}\]-1\)+\[{r\over2}\]_c+1\)+\(\[{t\over2}\]_c+1\)
=
\cr
&=
\(r\({N-1\over2}-1\)+\[{r\over2}\]_c+1\)+\(\[{n-rN\over2}\]_c+1\)
=
\cr
&=
r\({N-1\over2}-1\)+{r\over2}+1+{n-rN\over2}+1
-\Delta_0(r,n-rN)-{1\over2}\Delta_{odd}(r,n-rN)
=
\cr
&=
{n\over2}-r+2-\Delta_0(r,n-rN)-{1\over2}\Delta_{odd}(r,n-rN)
=
\cr
&=
\[{n\over2}\]-r+2-\Delta_0(r,n-rN)-
{1\over2}\Bigl(\Delta_{odd}(r,n-rN)
-
\Delta_{odd}(n)\Bigr)
=
\cr
&=
\[{n\over2}\]-r+2-\Delta_0(r,n-rN)-
\Delta_{odd}(r)\Bigl(1-\Delta_{odd}(n)\Bigr)
=
\cases{
\[{n\over2}\]-r+1&if a), b), or c) holds;
\cr\cr
\[{n\over2}\]-r+2&otherwise.
}
}
$$
Comparing this with Corollary 1,
we conclude that Theorem 1 is proven (modulo the main theorem).

\s 2.
Howie diagrams

Suppose that $S$ is a closed oriented
surface (possibly non-connected), and
$\Gamma$ is a finite
(undirected)
graph embedded into $S$ and dividing it into
simply connected
domains.
Such a graph determines a cell decomposition of $S$,
i.e. a mapping $\rm M$
called a \emph{map} on $S$:
$$
{\rm M}\:\bigsqcup\limits_{i=1}^m D_i\to S,
\qbox{where $D_i$ are two-dimensional disks,}
$$
such, that
\-
the mapping
$\rm M$ is continuous surjective,
injective on the interior
(i.e. on $\bigsqcup\limits_{i=1}^m(D_i\setminus\d D_i)$);
\-
the preimage of each point is finite, and
the preimage of the graph $\Gamma$ is the union of the boundaries
of the faces:
${\rm M}^{-1}(\Gamma)=\bigsqcup\limits_{i=1}^m \d D_i$.

\enditem
The preimages of the vertices of $\Gamma$
are
called \emph{corners} of the map;
we say that a corner $c$ is
\emph{at} a vertex~$v$ if ${\rm M}(c)=v$.
The vertices and edges of $\Gamma$ are referred to as
\emph{vertices and edges} of the map $\rm M$.
The disks $D_i$
are called \emph{faces} or \emph{cells}
of the map.
Such a map is
called a \emph{diagram} over a free product $A*B$ if
\-
the graph $\Gamma$ is
bipartite, i.e. there are two types of vertices:
$A$-vertices and
$B$-vertices, and each edge joins an~$A$-vertex with a $B$-vertex;
\-
the
corners at $A$-vertices are labeled by elements of
the group $A$, and the corners
at
$B$-vertices are labeled by elements of $B$;
\-
some vertices are distinguished and called \emph{exterior},
the
other vertices
are called \emph{interior}\/;
\-
the label of each interior $A$-vertex equals 1 in the group $A$, and
the label of each interior $B$-vertex equals 1 in $B$, where the
\emph{label of a vertex} is the product of labels of corners at this
vertex in clockwise order (thus,
the label of a vertex is defined up to conjugation in $A$ or
$B$).

\enditem
Similar diagrams were considered in
[How83],
[How90],
[K93],
[Le09],
and many other works, but our definitions slightly differ
and corresponds to the definitions from [IK18]
(except that exterior and interior vertices
are called irregular and regular in [IK18]).

The \emph{label of a face} of a diagram
is the
product of labels of all corners of this face in counterclockwise
order.
The label of a face is an element of the free product $A*B$
defined up to conjugation.

For instance, Figure 1
shows a diagram on a torus (which is drawn as a rectangle with
identified opposite sides) containing two vertices, three
edges, one face, and six corners with labels $a\in A$ and $b\in B$.
If both vertices are interior, then $a^3$ must be equal to 1 in
$A$, and $b^3$ must be equal to 1 in $B$. The label
of the face is
$(ab)^3$.
Actually, this diagram shows that the cube of the product
of two elements of
order three is always a commutator (in any group).

\vfill\break

\goodbreak
\bigskip
\centerline{\input 1.PIC}
\nobreak%
\centerline{Fig. \lowercase{1}}%
\goodbreak
\bigskip

\s 3.
Motions

This section is very similar to corresponding sections of
[FK12] and [IK18] and contains
definitions and statements from
[Kl05] with some simplifications.

Let $\Mu$ be a map on a closed oriented surface $S$
and let $\Gamma\subset S$ be the corresponding graph.
A \emph{car} moving around a face~$D$ is an
orientation preserving homeomorphism
from an oriented circle~$R$ (the \emph{circle of time})
to the boundary~$\d D$ of $D$.

If the number of cars being at a moment of time $t$ at a point $p$
of $\Gamma$
equals the degree $d$ of this point, then we
say that
a \emph{complete collision} (of degree~$d$) occurs at
$p$ at the moment $t$;
this point $p$ is called a \emph{point of complete collision}.
Here, the \emph{degree} of a point $p\in\Gamma$
is the number of edges incident to $p$
if $p$ is a vertex;
and $\deg p\:=2$ if $p$ is not a vertex
(i.e. if $p$ is an interior points of an edge).

Note that, according to the definition,
when a car arrives to vertex of degree one
(a dead end), a complete collision occurs.

\emph{A multiple motion} of period $T$ on a map $\Mu$ is a tuple of
cars $\alpha_{D,j}\:R\to\d D$, where $j=1,\dots,d_D$, such that

\item{1)}
$d_D\ge 1$ for any face $D$ (i.e. each face is moved around by at least one
car);

\item{2)}
$\alpha_{D,j}(t+T)=\alpha_{D,j+1}(t)$ for any $t\in R$ and
$j=\{1,\dots,d_D\}$ (here indices are modulo ${d_D}$,
and the addition of points of the circle $R$ is defined naturally:
$R=\R/P\Z$);

\item{3)}
for every face $D$, there exists a partition of $\d D$ into $d_D$
consecutive arcs with disjoint interiors such that, during the time
interval $0\le t\le T$, each car $\alpha_{D,j}$ is moving along the
$j$th arc of the partition.

\enditem
Informally,
several ($d_D$) cars are moved around
each face $D$
in
counterclockwise
direction
(the interior of $D$ remains on the left)
without U-turns and stops; and the motion is periodic in the sense
that the boundary of $D$ is partitioned into $d_D$ segments,
and, during the
period (of $T$ minutes), each car is moving along its segment (thus,
after $T$ minutes,
the cars' positions interchange cyclically).

\proclaim Car-crash lemma
{\rm(for multiple motions) [Kl05], [Kl97]}.
For any multiple motion
on a map on a closed oriented
surface $S$, the
number of points of complete collision is at least
$
\chi(S)+\sum\limits_D(d_D-1),
$
{where
the summation runs over all
faces of the map.}

In [Kl05] and [Kl97],
this lemma was stated and proven for connected surfaces,
but it remains valid in non-connected case because
the both sides of the inequality are additive with respect to
the disjoint union.

Consider, for instance, the following motion on the one-cell map on
a
torus shown in Figure 1:
three cars move around the unique face
with constant speed one edge per minute; at zero moment
of time, these
three cars are at three different angles with label $a$.
Figure 1
shows the location of cars at the moment $t=1/3$. This
is a periodic motion with period two minutes. Complete collisions
occur at
both vertices; while outside the vertices (i.e. at interior points
of edges) there are no
collisions.
The car-crash
lemma says that the following inequality must hold:
$$
\small
\pmatrix{
\hbox{the number of points of complete collision,}
\cr
\hbox{i.e. 2}
}
\ge
\pmatrix{
\hbox{the Euler characteristic of the torus,}
\cr
\hbox{i.e. 0}
}
+
\pmatrix{
\hbox{$d_D$, i.e. the number of cars}
\cr
\hbox{moving around the unique face $D$,}
\cr
\hbox{i.e. 3}
}
-1,
$$
which appears to be an equality in this example.


\s 4. Clusters

The idea of clusters is that collisions
that occur near each other can be treated as
one collision; the modified car-crash lemma
(the cluster lemma below) says that not only
the number of points of collision is large,
but also the number of points of collision that are far from each other is
large.

Suppose that we have a multiple motion with period $T$
on some map
on a surface $S$
and all cars move with the same constant speed one edge per minute.
A set $K$ of points of complete collision is called a \emph{cluster
centred at $v\in K$} if,
during less than $T/2$ minutes
after the collision at $v$,
each point
$w\in K$
is
visited by at least one car
having collided at $v$.
The cars
colliding at the centre $v$ of a cluster $K$
are referred to as the
\emph{connecting cars} of $K$;
\emph{the connecting paths} of $K$
are
the paths (of length $<T/2$) the connecting cars
move along on the way from the centre of $K$
to other points of the cluster.
A set $\cal C$ of clusters are called \emph{independent}
if the centre of each cluster from $\cal C$
does not lie on any connecting path of another cluster from
$\cal C$.

The statement of the cluster lemma (see below) uses the
\emph{fair partition function} $\fp(\M)$ of a multiset~$\M$
consisting of positive integers:
$$
\fp(\M)\:=\min\left\{\max\(\sum_{i\in\A}i,\;
\sum_{i\in\M\setminus\A}\!\!\!i\)
\;\Biggm|\;\A\subseteq\M\right\}.
\qbox{
For example, $\fp(10,4,4,3,2)=\max(10+2,\;4+4+3)=12$.
}
$$
The problem of finding a fair partition is
sometimes called ``the easiest
NP-hard problem"\ [Me06].
We need a simple example of such calculation:
$$
\fp(\,\overbrace{\underbrace{1,1,1,1,1\dots,1}_
{\min(l,\kappa)\hbox{ \small ones}},
N,N,\dots,N}^{\kappa\hbox{ \small numbers}}\,)
=
\cases{
\[{\kappa+1\over2}\],&if $\kappa\le l$;
\cr\cr
\[{l+1\over2}\]+N\cdot{\kappa-l\over2},
&if $\kappa>l$ and $\kappa-l$ is even;
\cr\cr
\[{l+1-\min(l,N)\over2}\]+N\cdot{\kappa-l+1\over2},
&if $\kappa>l$ and $\kappa-l$ is odd.
\cr
}
\eqno{(1)}
$$
This is true for all  $\kappa,N\in\N$ and $l\in\N\cup\0$.
The point is that
the following algorithm gives a fair partition $\M=\A\sqcup\B$:
\-
divide large items (i.e. $N$s) fairly, i.e.
give
$[(\kappa-l)/2]$ of these items to~$\A$
(if $\kappa>l$);
\-
use small items (i.e. ones) to
compensate the difference between $\A$ and $\B$ (which arises for
odd $\kappa-l$);
\-
when (and if) the difference is compensated, divide the remaining 1s
fairly.

\enditem
We leave
the proof
to the reader as an easy exercise;
Figures 2 and 3 show all possible cases
\newline
($f$ denotes
$\fp(\,\overbrace{\underbrace{1,1,1,1,1\dots,1}_
{\min(l,\kappa)\hbox{ \small things}},
N,N,\dots,N}^{\kappa\hbox{ \small things}}\,)
$
in these figures).

\goodbreak
\bigskip
\centerline{\input 2.PIC}
\nobreak%
\centerline{Fig. \lowercase{2}}%
\goodbreak
\bigskip

\vfil\break

\goodbreak
\bigskip
\centerline{\input 3.PIC}
\nobreak%
\centerline{Fig. \lowercase{3}}%
\goodbreak
\bigskip

We call a multiple motion \emph{uniform} if all cars
move with same constant speed one edge per minute
and are at some vertices at the moment $t=0$.

\proclaim Cluster lemma.
Suppose that,
for a multiple
uniform
motion
on a map
on a closed oriented surface $S$,
the set $\Pi$ of points of complete collision
is partitioned into the minimal possible
number $\kappa$ of
independent
clusters\: $\Pi=\bigsqcup\limits_{i=1}^\kappa K_i$.
Then
\item{\rm a)}
$
\kappa\ge
\chi(S)+\sum\limits_D(d_D-1),
$
where
the summation runs over all
faces\;
\item{\rm b)}
if
there are precisely $n$ points of
complete collision, and their degrees are
$N_1\le\dots\le N_n$,
then the number of cars of this motion \(i.e.
$\sum\limits_D d_D$,
where
the sum runs over all
faces\)
satisfies the inequality
$
\sum\limits_D d_D\ge
\max\bigl(
\fp(N_1,\dots,N_\kappa),
\;
N_n
\bigr);
$
in particular, if all points of complete collision have
degree at least $N$, then
$
\sum\limits_D d_D
\ge
\[{\kappa+1\over2}\]\cdot N.
$

\Proof
We assume that all collisions occur at vertices. This can be achieved by
the subdivision of each edge of the initial map into two equal parts by
new vertices of degree two (and slowing down all cars).

Let us prove the first assertion. For each cluster $K=K_i$ centred at
$v=v_i$, consider a minimal set of connecting paths $\pi_j=\pi_{i\,j}$
such that these paths contain all points of $K$. Due to minimality, for
each connecting car, we have at most one corresponding connecting path
$\pi_j$ lying in the boundary of a cell $D_j=D_{i\,j}$, which is moved
around by this car. By the definition of cluster, the length
$\tau_j$ of the path $\pi_j$
is less than $T/2$.

Let us connect the starting and the ending points of the path $\pi_j$
by a new path $\pi_j'$ of the same length lying inside
the  cell $D_j$ (so, we duplicate the path $\pi_j$).
(Note that these dublications for all clusters under
consideration may produce several chords
inside a cell, but these chords never intersect,
because the clusters are independent.)

The cell
$D_j$ turns into two cells (see Figure 4, on the left):
the \emph{large} cell $D_j'$ of the same perimeter as the initial cell
$D_j$
and the \emph{small} cell $\Gamma_j$ of perimeter
$2\tau_j$.


\goodbreak
\bigskip
\centerline{\input 4.PIC}
\nobreak%
\centerline{Fig. \lowercase{4}}%
\goodbreak
\bigskip

We want to define
car motion on this modified map.
Cars
moving around the
large cells imitate
the cars
moving around the
initial cells $D_j$, except that
they use the new road $\pi_j'$
instead of $\pi_j$.

To define a car moving around
the small cell, let us
consider the motion
of
already defined cars on the boundary of a small cell.
The boundary of each small cell~$\Gamma=\Gamma_j$
has length
$2\tau_j$
and
consists of three segments
(listed counterclockwise, see Figure~4, on the right):

\-
{part ${\bf a}=\pi_j'$}
of length
$\tau=\tau_j$;
in this segment the connecting car is moving
during
time interval $0\le t\le\tau$
(to simplify notation we assume that the complete collision
at $v_i$ occurs at zero moment, other cases
can of course be considered similarly);
part $\bf a$
ends at the corner $c$ at the centre $v_i$ of the cluster~$K_i$;

\-
{part $\bf b$} of
length one (this is the first edge of the path $\pi_j$),
starting at the corner $c$;
in this segment, another connecting car
of the cluster $K_i$ is moving
during
time $-1\le t\le 0$;

\-
{part $\bf c$}
of length $|{\bf c}|=\tau-1$;
not much is known about
cars moving here;
however, we know that
$\tau<T/2$, therefore,
$|{\bf c}|=\tau-1<T-\tau-1$;
this means that

\disp{\sl
there are no cars in
$\bf c$ \(including its ends\)
at some moment of time
$\tau<t<T-1$
and even during some subinterval $\Delta_\Gamma$
\(of positive duration\) of the time interval $\tau<t<T-1$
}
\item{}
(because the
time interval $\tau<t<T-1$ has
duration
$T-1-\tau>|{\bf c}|$,
and all cars move with the unit speed).
Note also that we can choose time intervals
$\Delta_\Gamma$ disjoint
for different small cells $\Gamma$:
$$
\Delta_\Gamma\cap \Delta_{\Gamma'}=\emptyset
\qbox{for }\Gamma\ne\Gamma'.
$$

\noindent
Now, let us define a new car
$\alpha_\Gamma$ moving around the small cell $\Gamma=\Gamma_j$:
\-
at moment zero, $\alpha_\Gamma$ is at the corner $c$
(and participates in the complete collision at $v_i$);
\-
then $\alpha_\Gamma$ moves (slowly) along the segment $\bf b$
without any  collisions, because the
(connecting) car moving in the opposite direction
has left segment
$\bf b$
having met our car $\alpha_\Gamma$ at $v_i$;
so segment $\bf b$ is safe until the moment $T-1$;
\-
during the time interval $\Delta_\Gamma$
(which starts earlier than $T-1$
by definition of $\Delta_\Gamma$),
our new car~$\alpha_\Gamma$ (rapidly)
moves through the segment $\bf c$;
no collisions occur, because $\bf c$ is safe
during the time interval $\Delta_\Gamma$ by definition of $\Delta_\Gamma$
(and because $\Delta_\Gamma\cap \Delta_{\Gamma'}=\emptyset$
for $\Gamma\ne\Gamma'$);
\-
thus, our car $\alpha_\Gamma$ arrives to segment
$\bf a$ later than moment $\tau$ (again, by definition
of $\Delta_\Gamma$);
this
means that the connecting car
already left $\bf a$, and our car $\alpha_\Gamma$ safely
without any collisions arrives to corner $c$ at the end of period.

\enditem
We obtain a periodic motion on
a map on surface $S$,
the number of
complete collisions is precisely $\kappa$,
and the sum~$\sum\limits_D(d_D-1)$ over all faces
remains the same as for the initial map
(because each small face $\Gamma_j$ is moved around by one car,
i.e. $d_{\Gamma_j}=1$).
Thus, applying the car-crash lemma to
this motion, we obtain assertion a).

\smallskip


\medskip
To prove assertion b),
we divide the period of time $I=\{t\;|\;0\le t<T\}$
into two half-periods: $I=I_1\sqcup I_2$, where
$I_1=\{t\;|\;0\le t<T/2\}$ and~$I_2=\{t\;|\;T/2\le t<T\}$.

The periodicity of the motion
implies that, at each point,
not more than one complete collision occurs
during the period $I$.
Therefore,
the set of points of complete collision
$\Pi=\{p_1,\dots,p_n\}$
is partitioned into two subsets:
$\Pi=\Pi_1\sqcup \Pi_2$, and
the multiset $\NN=(N_1,\dots,N_n)$ of the degrees
of these points
is partitioned into two submultisets:
$$
\NN=\NN_1\sqcup\NN_2,
\qbox{where
$\NN_i=\(\hbox{degrees of points of complete collision
occurring during $I_i$}\)$}.
$$

Suppose that $\Pi_i$ can be partitioned into
$\kappa_i$
independent
clusters and cannot be partitioned into
a fewer number of
independent
clusters.
Then $\kappa_1+\kappa_2\ge\kappa$
(because $\Pi$ can be partitioned into $\kappa_1+\kappa_2$
independent
clusters,
which is impossible if $\kappa_1+\kappa_2<\kappa$).

We say that a set of points of complete collision
is \emph{independent} if the sets of colliding cars at these points
during the period $I$ are disjoint.

\noindent
Let us concentrate on $\Pi_1$ now. Suppose that
\-
$v_1\in\Pi_1$ is a point at which the first
(timewise)
collision occurs (if such $v_1$ exists);
\-
$v_2\in\Pi_1$ is a point at which the first (timewise)
collision occurs such that the set $\{v_1, v_2\}$
is
independent
(if such~$v_2$ exists);
\-
$v_3\in\Pi_1$ is a point at which the first (timewise)
collision occurs such that the set $\{v_1, v_2, v_3\}$
is independent
(if such $v_3$ exists);
\-
\dots

\noindent
The number of points $v_i$ is at least $\kappa_1$, because
otherwise $\Pi_1$ would admit a partition on a less than
$\kappa_1$ number of
independent
clusters
(e.g., if $v_1$ and $v_2$ do exist, and $v_3$ does not,
then each point from $\Pi_1$ is either in a cluster centred at $v_1$
or in a cluster centred at $v_2$).

Thus, $\Pi_1$ contains independent points
$v_1,\dots,v_{\kappa_1}$, and $\Pi_2$
contains independent points $w_1,\dots,w_{\kappa_2}$
(by similar reasons).
Therefore, the number of all existing cars
is at least
$$
\max\(\sum\deg v_i,\;\sum\deg w_i\)
\ge
\fp(\deg v_1,\dots,\deg v_{\kappa_1}, \deg w_1,\dots,\deg w_{\kappa_2})
\ge
\fp(N_1,\dots,N_\kappa)
$$
(where the last estimate follows immediately from the inequalities
$\kappa_1+\kappa_2\ge\kappa$ and $N_1\le\dots\le N_n$).

This completes the proof of assertion b), because the bound
$\sum\limits_D d_D\ge N_n$ is obvious
(if, at some point, $N_n$ cars collide,
then $N_n$ cars do exist).

\medskip\noindent
The cluster lemma implies the following fact
(not mentioning clusters at all).

\Corollary of the cluster lemma.
Suppose that a multiple
uniform
motion
on a map
on an oriented closed
surface~$S$
has precisely $n$ points of
complete collision, and their degrees are
$N_1\le\dots\le N_n$.
Then the number of cars of this motion \(i.e.
$\sum\limits_D d_D$,
where
the summation runs over all
faces\)
satisfies the inequality
$$
\sum\limits_D d_D\ge
\max\bigl(
\fp(N_1,\dots,N_\kappa),
\;
N_n
\bigr),
\qbox{where
$\kappa=\chi(S)+\sum\limits_D(d_D-1)$
{\rm(this value never exceeds $n$).}}
\eqno{(2)}
$$
Moreover, for all $l\in\N\cup\0$,
$$
\chi(S)-l+\sum_D(d_D-1)
\le
\cases{
2\[\frac1{N_{l+1}}\(\sum\limits_D d_D-\[\frac{l+1}2\]\)\],
&if $\sum\limits_D(d_D-1)-l$ is even\;
\cr\cr
2\[\frac1{N_{l+1}}\(\sum\limits_D d_D-\[\frac{l+1-N_{l+1}}2\]_+\)\]-1,
&if $\sum\limits_D(d_D-1)-l$ is odd,
}
\eqno{(3)}
$$
where $[x]_+\:=\max([x],0)$ and $N_i=\infty$ for $i>n$
\(in particular, for $N_{l+1}=\infty$,
the right-hand side is $0$ or $-1$\).

\Proof
To prove (2),
it suffice to
substitute bound a) of the
cluster lemma to
bound b) of the same lemma
(as
the fair partition function $\fp(N_1,\dots,N_\kappa)$
is surely non-decreasing as a function of $\kappa$).

Let us prove (3).
Using the monotonicity of the function $\fp$
with respect to each argument and
formulae (2)~and~(1),
for
$\kappa=\chi(S)+\sum\limits_D(d_D-1)$,
we obtain
$$
\eqalign{
\sum\limits_D d_D
\gee^2&
\max\bigl(
\fp(N_1,\dots,N_\kappa),
\;
N_n
\bigr)
\ge
\fp(\,\overbrace{\underbrace{1,1,1,1,1\dots,1}_
{\min(l,\kappa)\hbox{ \small ones}},
N_{l+1},N_{l+1},\dots,N_{l+1}}^{\kappa\hbox{ \small numbers}}\,)
\=^1
\cr\cr
\=^1&
\cases{
\[{\kappa+1\over2}\]&if $\kappa\le l$;
\cr\cr
\[{l+1\over2}\]+N_{l+1}\cdot{\kappa-l\over2}
&if $\kappa>l$ and $\kappa-l\in2\Z$;
\cr\cr
\[{l+1-\min(l,N_{l+1})\over2}\]+N_{l+1}\cdot{\kappa-l+1\over2}
&if $\kappa>l$ and $\kappa-l\notin2\Z$.
\cr
}
}
$$

\Case 0: $\kappa\le l$.
\-
If $\kappa-l$ is even, then
$
[(\kappa+1)/2]
\ge
[(l+1)/2]+N_{l+1}\cdot(\kappa-l)/2
$
(since $N_{l+1}\ge1$).
\-
If $\kappa-l$ is odd, then
$
[(\kappa+1)/2]
\ge
[l/2]+N_{l+1}\cdot(\kappa-l+1)/2
\ge
[(l+1-\min(l,N_{l+1}))/2]+N_{l+1}\cdot(\kappa-l+1)/2
$.
\enditem
\medskip\noindent
Thus, for all $\kappa$ and $l$,
$$
\sum\limits_D d_D
\ge
\cases{
\[{l+1\over2}\]+N_{l+1}\cdot{\kappa-l\over2}
&if $\kappa-l\in2\Z$;
\cr\cr
\[{l+1-\min(l,N_{l+1})\over2}\]+N_{l+1}\cdot{\kappa-l+1\over2}
&if $\kappa-l\notin2\Z$.
\cr
}
$$

\Case 1: $\kappa-l$ is even.
$$
\sum\limits_D d_D
\ge
\[{l+1\over2}\]+N_{l+1}\cdot{\kappa-l\over2}
\ \imp\
\kappa-l
\le
{2\over N_{l+1}}\(\sum\limits_D d_D -\[{l+1\over2}\]\)
\ \imp\
\kappa-l
\le
2\[{1\over N_{l+1}}\(\sum\limits_D d_D -\[{l+1\over2}\]\)\],
$$
where the last implication is valid since $\kappa-l\in2\Z$.
The obtained bound coincides with (3), because
$\kappa=\chi(S)+\sum\limits_D(d_D-1)$.

\Case 2: $\kappa-l$ is odd.
$$
\sum\limits_D d_D
\ge
\[{l+1-N_{l+1}\over2}\]_++N_{l+1}\cdot{\kappa-l+1\over2}
\ \imp\
\kappa-l+1
\le
{2\over N_{l+1}}\(\sum\limits_D d_D -\[{l+1-N_{l+1}\over2}\]_+\);
$$
$\kappa-l+1$
is even, hence,
$
\kappa-l+1
\le
2\[{1\over N_{l+1}}\(\sum\limits_D d_D -\[{l+1-N_{l+1}\over2}\]_+\)\],
$
as required. This completes the proof.

\s 5.
Main theorem

\proclaim Main theorem.
Suppose that, in a
free product of groups $G=\zvezda\limits_{j\in J} A_j$,
an equality
$$
c_1\dots c_kd_1\dots d_l=u_1^{n_1}\dots u_m^{n_m}
\eqno{(4)}
$$
holds,
where $c_i$ are commutators, $d_i$ are conjugate to elements of
$\bigcup\limits_{j\in J}A_j$,
elements
$u_i$ are conjugate to each other and not conjugate to elements of
$\bigcup\limits_{j\in J}A_j$,
and $n_i$ are positive integers. Then
$$
2-2k-l+\sum_{i=1}^m(n_i-1)
\le
\cases{
2\[\frac1N\(\sum\limits_{i=1}^m n_i -\[\frac{l+1}2\]\)\],
&if $\sum\limits_{i=1}^m(n_i-1)-l$ is even\;
\cr\cr
2\[\frac1N\(\sum\limits_{i=1}^m n_i -\[\frac{l+1-N}2\]_+\)\]-1,
&if $\sum\limits_{i=1}^m(n_i-1)-l$ is odd,
}
$$
where $[x]_+\:=\max([x],0)$, and
$N$ is the minimal order of
a letter \(from $\bigcup\limits_{j\in J}A_j$\) of a
cyclically
reduced word $u$ conjugate to $u_1$
\(in particular, for $N=\infty$,
the right-hand side is $0$ or $-1$\).

\Proof
Without loss of generality, we can assume that the
number $|J|$ of free factors is two.
Indeed, suppose that the cyclically irreducible form $u$
of the elements $u_i$ contains a letter $a_j\in A_j\setminus\1$
for some $j\in J$. Then $G$ decomposes into the free
product $G=A*B$, where $A=A_j$ and
$B=\zvezda\limits_{j'\in J\setminus\{j\}} A_{j'}$,
and the conditions of the theorem remain fulfilled for this decomposition.
Thus, we assume that the number of factors is two: $G=A*B$.

Equality (4) allows us to draw a Howie diagram
on an oriented closed
(possibly non-connected)
surface
$S$ of genus~$k'\:={1\over2}(2-\chi(S))$
with~$l'$ exterior vertices and $m$ cells, whose labels are
$u^{n_1},\dots,u^{n_m}$,
where
$$
k'\le k
\qqbox{and}
2k'+l'\le2k+l.
$$
In detail, this construction
was explained in [IK18]. Let us restrict ourselves to
just one example. If $k=m=1$ and $l=0$, then, in most cases,
we naturally obtain
a torus without exterior vertices (with several interior vertices), and
with one cell, whose label is $u_1^{n_1}$ (see Figure~1, for instance);
but if, e.g., commutator $c_1$ has the form $c_1=[a,v]$, where
$a\in A\setminus\1$ and $v\in(A*B)\setminus A$, then we obtain
a sphere with two
exterior vertices (whose labels are $a$ and $a^{-1}$) and one cell, whose
label is $c_1$. For $m>1$, we can even obtain a non-connected surface
(if equality~(4) decomposes into a product of two equalities of the same
type).

On the obtained diagram,
a multiple uniform motion is
naturally defined: a cell with label $u^{n_i}$ are moved around by $n_i$
cars with unit speed (one edge per minute); at moment $s\in\Z$,
each car is at the corner whose label is the $s$th letter of
the word $u$ (where $s$ is counted modulo the length of $u$).

\noindent
Thus,
collisions outside vertices (i.e. inside edges)
cannot occur, because at each moment of time either
\-
all cars are at $A$-vertices,
\-
or all cars are at $B$-vertices,
\-
or each car is moving along an edge from an $A$-vertex to a $B$-vertex,
\-
or each car is moving along an edge from a $B$-vertex to an $A$-vertex.

\enditem
A complete collision at a vertex
$v$ means that all corners at this vertex have the same label equal to
a letter of the word $u$. If vertex $v$ is interior, then
the product of these labels must be 1, i.e. $\deg v\ge N$.
Applying the corollary of the cluster lemma to this motion, we obtain
the inequality
$$
\Phi(k',l')\:=2-2k'-l'+\sum_{i=1}^m(n_i-1)
\le
\Psi(l')\:=
\cases{
2\[\frac1N\(\sum\limits_{i=1}^m n_i -\[\frac{l'+1}2\]\)\],
&if $\sum\limits_{i=1}^m(n_i-1)-l'\in2\Z$\;
\cr\cr
2\[\frac1N\(\sum\limits_{i=1}^m n_i -\[\frac{l'+1-N}2\]_+\)\]-1,
&if $\sum\limits_{i=1}^m(n_i-1)-l'\notin2\Z$.
}
\eqno{(5)}
$$
Note that
$$
\Psi(l+2)\le\Psi(l)
\qqbox{and}
\Psi(l\pm1)\le\Psi(l)+1
\qbox{for all $l$}.
\eqno{(6)}
$$
The first inequality is obvious;
to explain the second one, we
put $n\:=\sum n_i$. Now, if~$\sum(n_i-1)-l$ is even, then
$$
\eqalign{
&\Psi(l\pm1)
\le
\Psi(l-1)=
2\[\frac1N\(n-\[\frac{l-N}2\]_+\)\]-1
\le
2\[\frac1N\(n-\[\frac{l-N}2\]\)\]-1
\le
\cr
&\le
2\[\frac1N\(n-\[\frac{l+1-2N}2\]\)\]-1
=
2\[\frac1N\(n-\[\frac{l+1}2\]+N\)\]-1
=
2\[\frac1N\(n-\[\frac{l+1}2\]\)\]+2-1
=
\Psi(l)+1;
}
$$
if $\sum(n_i-1)-l$ is odd, then
$
\Psi(l\pm1)
\le
\Psi(l-1)
=
2\[\frac1N\(n-\[\frac{l}2\]\)\]
\le
2\[\frac1N\(n-\[\frac{l+1-N}2\]_+\)\]
=\Psi(l)+1.
$
This proves~(6).

Now,
if $l'\ge l$, then
$\Phi(k,l)\le\Phi(k',l')$, because $2k'+l'\le2k+l$
(as was noted above), and
$\Psi(l')\lee^6\Psi(l)+1$. Therefore,
$\Phi(k,l)\le\Phi(k',l')\le\Psi(l')\lee^6\Psi(l)+1$
and, hence, $\Phi(k,l)\le\Psi(l)$, because
$\Phi(k,l)$ and $\Psi(l)$ have the same parity (see (5)).
This completes the proof in the case where $l'\ge l$.

Now, suppose that $l'<l$. Note that~(6) implies
that
the function $l\mapsto l+\Psi(l)$ is non-decreasing. Thus,
the assertion of the theorem follows immediately
from~(5), because $k'\le k$ (as is noted above).

\baselineskip 10.5pt

\References

[Cal09]
D. Calegari,
scl.
MSJ Memoirs, 20.
Mathematical Society of Japan, Tokyo, 2009. xii+209 pp.

[Ch18]
L. Chen,
Spectral gap of scl in free products,
Proc. Amer. Math. Soc., 146:7 (2018), 3143-3151.
\arXiv 1611.07936

[Cl03]
A. Clifford,
Non-amenable type K equations over groups,
Glasgow Mathematical Journal 45:2 (2003), 389-400.

[ClG95]
A. Clifford and R. Z. Goldstein,
Tesselations of $S^2$ and equations over torsion-free groups,
Proceedings of the Edinburgh Mathematical Society 38:3 (1995), 485-493.

[ClG01]
A. Clifford and R. Z. Goldstein,
The group $\pres<G,t|e>$ when $G$ is torsion free,
Journal of Algebra 245:1 (2001), 297-309.

[CoR01]
M. M. Cohen and C. Rourke,
The surjectivity problem for one-generator,
one-relator extensions of torsion-free groups,
Geometry \& Topology 5:1 (2001), 127-142.
\arXiv math/0009101

[CCE91]
J. A. Comerford, L. P. Comerford Jr. and C. C. Edmunds,
{Powers as products of commutators},
Comm. Algebra, {19:2} (1991), 675-684.

[CER94]
L. P. Comerford Jr., C. C. Edmunds and G. Rosenberger,
Commutators as powers in free products of groups,
Proc. Amer. Math. Soc., 122:1 (1994), 47-52.
\arXiv math/9310205


[Cull81]
M. Culler,
Using surfaces to solve equations in free groups,
Topology, 20:2 (1981), 133-145.

[DH91]
A. J. Duncan and J. Howie,
{The genus problem for one-relator products of locally indicable
groups},
Mathematische Zeitschrift, {208:1} (1991), 225-237.

[FeR96]
R. Fenn and C. Rourke,
Klyachko's methods and the solution of equations over torsion-free groups,
\newline
L'Enseignement Math\'ematique 42 (1996), 49-74.

[FoR05a]
M. Forester and C. Rourke,
Diagrams and the second homotopy group,
Communications in Analysis and Geometry 13:4 (2005), 801-820.
\arXiv math/0306088

[FoR05b]
M. Forester and C. Rourke,
The adjunction problem over torsion-free groups,
Proceedings of the National Academy of Sciences of
the USA 102:36 (2005), 12670-12671.
\arXiv math/0412274

[FK12]
E. V. Frenkel and Ant. A. Klyachko,
Commutators cannot be proper powers in metric
small-cancellation torsion-free groups,
arXiv:1210.7908.

[How83]
J. Howie,
The solution of length three equations over groups,
{Proc. Edinburgh Math. Soc.}, {26:1} (1983), 89-96.

[How90]
J. Howie,
The quotient of a free product of groups by a single high-powered relator.
II. Fourth powers,
Proc. London Math. Soc., s3-61:1 (1990), 33-62.

[IK00]
S. V. Ivanov and Ant. A. Klyachko,
Solving equations of length at most six over torsion-free groups,
Journal of Group Theory 3:3 (2000), 329-337.

[IK18]
S. V. Ivanov and Ant. A. Klyachko,
Quasiperiodic and mixed commutator factorizations in free products of groups,
Bull. London Math. Soc., 50:5 (2018), 832-844.
\arXiv 1702.01379

[K93]
Ant. A. Klyachko,
A funny property of a sphere and equations over groups,
{Comm. Algebra}, {21:7} (1993), 2555-2575.

[Kl97]
Ant. A. Klyachko,
Asphericity tests,
Internat. J. Algebra Comp., 7:4 (1997), 415-431.

[Kl05]
Ant. A. Klyachko,
The Kervaire--Laudenbach conjecture and presentations of simple groups,
{Algebra and Logic}, 44:4 (2005), 219-242.
\arXiv:math.GR/0409146

[Kl06a]
Ant. A. Klyachko,
How to generalize known results on equations over groups,
Mathematical Notes 79:3-4 (2006), 377-386.
\arXiv math.GR/0406382

[Kl06b]
Ant. A. Klyachko,
SQ-universality of one-relator relative presentations,
Sbornik: Mathematics 197:10 (2006), 1489-1508.
\arXiv math.GR/0603468

[Kl07]
Ant. A. Klyachko,
Free subgroups of one-relator relative presentations,
Algebra and Logic 46:3 (2007), 158-162.
\arXiv math.GR/0510582

[Kl09]
Ant. A. Klyachko,
The structure of one-relator relative presentations and their centres,
Journal of Group Theory, 12:6 (2009), 923-947.
\arXiv math.GR/0701308

[KlL12]
Ant. A. Klyachko and D. E. Lurye,
Relative hyperbolicity and similar
properties of one-generator one-relator relative presentations with
powered unimodular relator,
J. Pure Appl. Algebra, 216:3 (2012), 524-534.
\arXiv 1010.4220

[Le09]
Le Thi Giang,
The relative hyperbolicity of one-relator relative presentations,
Journal of Group Theory, 12:6 (2009), 949-959.
\arXiv 0807.2487

[Me06]
S. Mertens,
The easiest hard problem: number partitioning,
Computational Complexity and Statistical Physics, 125:2 (2006),
125-139.
\arXiv cond-mat/0310317

[Sch59]
M. P. Sch\"utzenberger, Sur l'equation $a^{2+n}=b^{2+m}c^{2+p}$ dans
un groupe libre,
C. R. Acad. Sci.  Paris S\'er. I Math., 248 (1959), 2435-2436.

\end